\documentclass[11pt]{article}
\usepackage{amssymb,amsmath,amsthm,amsfonts,amsxtra}
\usepackage{mathrsfs,amscd}
\makeatletter
\def\@fnsymbol#1{\ensuremath{\ifcase#1\or 1\or 2\fi}}
\makeatother

\def\shd{\mathcal{D}}

\def\shm{\mathcal{M}}
\def\sho{\mathcal{O}}

\newcommand{\C}{\mathbb{C}}
\newcommand{\N}{\mathbb{N}}
\newcommand{\R}{\mathbb{R}}

\newcommand{\Z}{\mathbb{Z}}
\newcommand{\K}{\mathbf{k}}
\newcommand{\I}{\mathcal{I}}
\newtheorem{theorem}{Theorem}[section]
\newtheorem{proposition}[theorem]{Proposition}
\newtheorem{lemma}[theorem]{Lemma}
\newtheorem{corollary}[theorem]{Corollary}

\theoremstyle{definition}

\title{Functorial properties of the microsupport and regularity for ind-sheaves}
\makeindex
\begin{document}

\author{Ana Rita Martins\footnote{The research of the author was supported by
Funda\c c{\~a}o para a Ci{\^e}ncia e Tecnologia and FEDER (project
POCTI-ISFL-1-143 of Centro de Algebra da Universidade de
Lisboa).}}

\date{\today}
\maketitle

\begin{abstract}
The notion of microsupport and regularity for ind-sheaves was
introduced by M. Kashiwara and P. Schapira in \cite{KS2}. In this
paper we study the behaviour of the microsupport under several
functorial ope-rations and characterize ``microlocally" the
ind-sheaves that are regular along involutive manifolds. As an
application we prove that if a cohe-rent $\shd$-module $\shm$ is
regular along an involutive manifold $V$ (in the sense of
\cite{KO}), then the ind-sheaf of temperate holomorphic solutions
of $\shm$ is regular along $V$. Another application is the notion
of microsupport for sheaves on the subanalytic site, which can be
given directly, without using the more complicated language of
ind-sheaves.
\end{abstract}

\tableofcontents

\section{Introduction}

\hspace*{\parindent}Let $X$ be a real manifold and $\K$ a field.
In \cite{KS3} M. Kashiwara and P. Schapira gave a new perspective
of sheaf theory with the introduction of the notion of
microsupport. The microsupport $SS(F)$ of a sheaf $F$ of
$\mathbf{k}$-vector spaces on a manifold $X$ is a closed conic
subset of the cotagent bundle $T^*X$ to $X$, describing the
directions of non propagation for $F$.

Later, motivated by the fact that various objects in Analysis
cannot be treated with sheaf theoretical methods, such as
functions with growth conditions, the category of sheaves was
enlarged to that of ind-sheaves and a corresponding notion of
microsupport was introduced in \cite{KS2}. Since this notion
coincides with the classical one for sheaves on $X$, it is natural
to ask if the functorial properties of the microsupport of
classical sheaves still hold for ind-sheaves. Our aim in this
paper is to study the behaviour of the microsupport under direct
image of closed embeddings, smooth inverse images and
$R\mathcal{I}{hom}(\cdot,\cdot)$. For
$R\mathcal{I}{hom}(\cdot,\cdot)$ we were only able to treat the
functor $R\mathcal{I}{hom}(q_2^{-1}(\cdot), q_1^{!}(\cdot))$,
where $q_1$ and $q_2$ denote the first and second projection
defined on $X\times Y$.

In \cite{KS2}, the authors gave an example showing that the
estimate $$SS(R\mathcal{H}{om}(G,F))\subset
SS(F)\widehat{+}SS(G)^a,$$ proved for classical sheaves, doesn't
work for ind-sheaves. This fact motivated the definition of
regularity for ind-sheaves. In this paper, we also study this
notion, characterizing ``microlocally" the ind-sheaves that are
re-gular along involutive manifolds. More precisely, we show the
following: Let $D^b(\text{I}(\K_X))$ be the bounded derived
category of $\text{I}(\K_X)$, the category of ind-sheaves on $X$,
$p$ an element of the cotagent bundle $\pi:T^*X\to X$ and
$D^b(\mathrm{I}(\K_X;p))$ the localization of
$D^b(\mathrm{I}(\K_X))$ by the subcategory of objects $F$ such
that $p\notin SS(F)$. Let $f:Y\to X$ be a morphism of manifolds.
When $f$ is a closed embedding and $p$ belongs to the conormal
bundle $T_Y^*X$ to $Y$ in $X$, then $F\in D^b(\text{I}(\K_X))$ is
regular along $\pi^{-1}(Y)$ at $p$ if and only if $F\simeq Rj_*G$
in $D^b(\mathrm{I}(\K_X;p))$, for some $G\in D^b(\text{I}(\K_Y))$.
When $f$ is smooth and $p\in Y\times_X T^*X$, identifying
$Y\times_X T^*X$ with a submanifold of $T^*Y$, then $G\in
D^b(\text{I}(\K_Y))$ is regular along $Y\times_X T^*X$ at $p$ if
and only if $G\simeq f^{-1} F$ in $ D^b(\mathrm{I}(\K_Y;p))$, for
some $F\in D^b(\text{I}(\K_X))$.

Let $X$ be a complex manifold and let $\sho_X^t$ denote the
ind-sheaf of temperate holomorphic functions on $X$. When $\shm$
is a coherent $\shd_X$-module regular along an involutive manifold
$V$, in the sense of \cite{KO}, we prove that the ind-sheaf
$Sol^t(\shm):=R\mathcal{I}{hom}_{\beta_X(\shd_X)}(\beta_X(\shm),
\sho_X^t)$ is regular along $V$. Moreover, when $F$ is an
$\R$-constructible sheaf, we may use the language of ind-sheaves
to give an alternative (and easy) proof of the estimate
$$SS(R\mathcal{H}{om}_{\shd_X}(\shm, {thom}(F,\sho_X))\subset V
\widehat{+}SS(F)^a,$$ given in \cite{KMS}.

Among all ind-sheaves, the ind-objects of the category of
$\R$-constructible sheaves are particularly interesting, since
they may be constructed using Grothendieck topologies. Let
$\text{Mod}_{\R-c}^c(\K_X)$ denote the category of
$\R$-cons-tructible sheaves with compact support. In \cite{KS1}
was given an equivalence between the category
$\text{I}\R-c(\K_X)=\text{Ind}(\text{Mod}_{\R-c}^c(\K_X))$ with
the category Mod$(\K_{X_{sa}})$ of sheaves on the subanalytic site
$X_{sa}$ associated to $X$. On the other hand, the natural exact
functor $\I:\text{I}\R-c(\K_X)\to \text{I}(\K_X)$ induces an
equivalence of triangulated categories
$\I:D^b(\text{I}\R-c(\K_X))\to
D_{\text{I}\R-c}^b(\text{I}(\K_X))$, where
$D_{\text{I}\R-c}^b(\text{I}(\K_X))$ denotes the full subcategory
of $D^b(\text{I}(\K_X))$ consisting of objects with cohomology in
I$\R-c(\K_X)$ (see \cite{KS1}). Therefore, identifying
$D^b(\K_{X_{sa}})$, the derived category of Mod$(\K_{X_{sa}})$,
with $D^b(\text{I}\R-c(\K_X))$, there is a natural way to
introduce the notion of microsupport for sheaves on $X_{sa}$,
setting $$SS(F)=SS(\I(F)),$$ for all
$F\in\text{Mod}(\K_{X_{sa}})$. In this paper we also translate
this definition of microsupport into several equivalent conditions
free of the so complicated language of ind-sheaves. We apply to
sheaves on $X_{sa}$ our results on the functorial properties of
the microsupport of ind-sheaves and prove the estimate:
$$SS(R\mathcal{H}{om}(q_2^{-1}G,q_1^{!}F))\subset SS(F)\times
SS(G)^a,$$ for $F\in D^b(\K_{X_{sa}})$ and $G\in
D^b(\K_{Y_{sa}})$, $q_1$ and $q_2$ denoting the first and second
projection defined on $X\times Y$.

We thank P. Schapira who suggested us to study these problems and
to T. Monteiro Fernandes for useful discussions through the
preparation of this work.

\section{Notations and review}

\hspace*{\parindent}We will mainly follow the notations in
\cite{KS1} and \cite{L}. In this paper, all manifolds will be real
analytic.\newline

\noindent\textbf{Geometry.} Let $X$ be a real analytic manifold.
We denote by $\pi:T^*X\rightarrow X$ the cotangent bundle to $X$
and we identify $X$ with the zero section of $T^*X$. Let $M$ be a
submanifold of $X$. We denote by $T_M X$ (resp. $T^*_M X$) the
normal bundle (resp. the conormal bundle) to $M$ in $X$.

Given a subset $A$ of $T^*X$, we denote by $A^a$ the image of $A$
by the antipodal map $$a:(x;\xi)\mapsto (x;-\xi).$$ The closure of
$A$ is denoted by $\overline{A}$. If $A$ is a locally closed
subset of $T^*X$, we say that $A$ is $\R^+$-conic (or simply
``conic" for short) if it is locally invariant under the action of
$\R^+$.

For a cone $\gamma\subset TX$, the polar cone $\gamma^\circ$ to
$\gamma$ is the convex cone in $T^*X$ defined by
\begin{center}
$\gamma^\circ=\{(x;\xi)\in T^*X; x\in\pi(\gamma)$ and $\langle
v,\xi\rangle\geq 0$ for any $(x;v)\in\gamma \}.$
\end{center}

When $X$ is an open subset of $\R^n$ and $\gamma$ is a closed
convex cone (with vertex at 0) in $\R^n$, we denote by $X_\gamma$
the open set $X$ endowed with the induced $\gamma$-topology of
$\R^n$ and by $\phi_\gamma$ the natural continuous map from $X$ to
$X_\gamma$ (see \cite{KS3}). In this paper, all closed convex
cones in $\R^n$ will be subanalytic.\newline

\noindent\textbf{Sheaves.} Let $X$ be a real analytic manifold and
$\K$ be a field. We denote by Mod$(\K_X)$ the abelian category of
sheaves of $\mathbf{k}$-vector spaces on $X$ and by
$D^b(\mathbf{k}_X)$ its bounded derived category. We denote by
Mod$_{\R-c}(\K_X)$ the abelian category of $\R$-constructible
sheaves of $\K$-vector spaces on $X$ and by $D^b_{\R-c}(\K_X)$ the
full subcategory of $D^b(\K_X)$ consisting of objects with
$\R$-constructible cohomology. Recall that the natural functor
$$D^b(\text{Mod}_{\R-c}(\K_X))\to D^b_{\R-c}(\K_X),$$ is an
equivalence of categories (see \cite{K1}). We denote by
Mod$_{\R-c}^c(\K_X)$ the full abelian subcategory of
Mod$_{\R-c}(\K_X)$ of $\R$-constructible sheaves with compact
support.

For an object $F\in D^b(\mathbf{k}_X)$, one denotes by $SS(F)$ its
\emph{microsupport}, a closed $\R^+$-conic involutive subset of
$T^*X$. We refer \cite{KS3} for details.\newline

\noindent\textbf{Ind-sheaves.} Let $X$ be a real analytic
manifold. One denotes by I$(\K_X)$ the abelian category of
ind-sheaves on $X$, that is, Ind(Mod$^c(\K_X)$), the category of
ind-objects of the category Mod$^c(\K_X)$ of sheaves with compact
support on $X$.

Recall the natural faithful exact functor
$$\begin{matrix}\iota_X:& \text{Mod}(\K_X) & \to &
\text{I}(\K_X)\\ & F& \mapsto &
\underset{\underset{U\subset\subset
X}{\longrightarrow}}{\sideset{``}{"}\lim }F_U \ (U \ \text{open}).
\end{matrix}$$ Thanks to $\iota_X$ we identify
Mod$(\K_X)$ with a full abelian subcategory of I$(\K_X)$ and
$D^b(\K_X)$ with a full subcategory of $D^b$(I$(\K_X)$).

The functor $\iota_X$ admits an exact left adjoint
$$\begin{matrix}\alpha_X:& \text{I}(\K_X)& \to &
\text{Mod}(\K_X)\\ & F=\underset{\underset{i\in I
}{\longrightarrow}}{\sideset{``}{"}\lim }F_i& \mapsto &
\underset{\underset{i\in I }{\longrightarrow}}{\lim }F_i.
\end{matrix}$$ This functor also
admits an exact fully faithful left adjoint
$\beta_X:\text{Mod}(\K_X)\to\text{I}(\K_X)$, which is not so easy
to describe. However, when $Z=U\cap S$, with $U$ open and $S$
closed subsets of $X$, one has $\beta_X(\K_Z)\simeq
\underset{\underset{S\subset W, V\subset\subset U
}{\longrightarrow}}{\sideset{``}{"}\lim }\K_{V\cap\overline{W}}$,
where $W$ (resp. $V$) runs through the family of open
neighborhoods of $Z$ (resp. relatively compact open subsets of
$U$).

We denote by $J$ the functor
$J:D^b(\text{I}(\K_X))\to\text{Ind}(D^b(\text{Mod}^c(\K_X)))$
defined by: $$J(F)(G)=\text{Hom}_{D^b(\text{I}(\K_X))}(G,F),$$ for
every $F\in D^b(\text{I}(\K_X))$ and $G\in D^b(\K_X)$.\newline

\noindent\textbf{Sheaves on the subanalytic site.} We denote by
$\text{Op}(X_{sa})$ the category of open subanalytic subsets of
$X$. One endows this category with a Grothendieck topology by
deciding that a family $\{U_i\}_i$ in $\text{Op}(X_{sa})$ is a
co-vering of $U\in \text{Op}(X_{sa})$ if, for any compact subset
$K$ of $X$, there exists a finite subfamily which covers $U\cap
K$. One denotes by $X_{sa}$ the site defined by this topology and,
for $U\in\text{Op}(X_{sa})$, we denote by $U_{sa}$ the category
$\text{Op}(X_{sa})\cap U$ with the topology induced by $X_{sa}$.
We have an equivalence of categories
$\text{Mod}(\K_{X_{sa}})\to\text{Mod}(\K_{X^c_{sa}})$, where
$\text{Mod}(\K_{X^c_{sa}})$ notes the category of sheaves on the
site $\text{Op}^c(X_{sa})$, consisting of relatively compact
subanalytic open subsets of $X$, with the topology induced by
$X_{sa}$.

The category of sheaves on $X_{sa}$ is denoted by
Mod$(\K_{X_{sa}})$. Recall that Mod$(\K_{X_{sa}})$ is a
Grothendieck category and, in particular, it has enough injective
objects.

Let $\rho:X\to X_{sa}$ be the natural morphism of sites. We have
functors
$$\text{Mod}(\K_X)\overset{\rho_*}{\underset{\rho^{-1}}{\rightleftarrows}}\text{Mod}(\K_{X_{sa}}),$$
and we still note by $\rho_*$ the restriction of $\rho_*$ to
Mod$_{\R-c}(\K_X)$ and to Mod$_{\R-c}^c(\K_X)$.

Let I$\R-c(\K_X)$ denote the category Ind(Mod$_{\R-c}^c(\K_X)$).
We may extend the functor
$\rho_*:\text{Mod}_{\R-c}^c(\K_X)\rightarrow\text{Mod}(\K_{X_{sa}})$
to $\text{I}\R-c(\K_X)$, by setting: $$\begin{matrix}\lambda: &
\text{I}\R-c(\K_X)&\to & \text{Mod}(\K_{X_{sa}})\\ &
\underset{\underset{i}{\longrightarrow}}{\sideset{``}{"}\lim }F_i
& \mapsto & \underset{\underset{i}{\longrightarrow}}{\lim
}\rho_*F_i,\end{matrix}$$  and it is proved in \cite{KS1} that
$\lambda$ is an equivalence of abelian categories. Most of the
time, thanks to $\lambda$, we identify $\text{I}\R-c(\K_X)$ with
$\text{Mod}(\K_{X_{sa}})$.

On the other hand, the natural functor
$\text{Mod}_{\R-c}^c(\K_X)\to \text{Mod}(\K_X)$ gives rise to a
functor $\I:\text{I}\R-c(\K_X)\to\text{I}(\K_X)$, which induces an
equivalence of triangulated categories:
\begin{equation}\label{SE:7}
\I:D^b(\text{I}\R-c(\K_X))\to D^b_{\text{I}\R-c}(\text{I}(\K_X)),
\end{equation}
where $D^b_{\text{I}\R-c}(\text{I}(\K_X))$ denotes the full
subcategory of $D^b(\text{I}(\K_X))$ consisting of objects with
cohomology in I$\R-c(\K_X)$.

Since the functor $\lambda$ induces an equivalence from
$D^b(\text{I}\R-c(\K_X))$ to the bounded derived category
$D^b(\K_{X_{sa}})$ of $\text{Mod}(\K_{X_{sa}})$, it gives rise to
an equi-valence still denoted by $\I:D(\K_{X_{sa}})\to
D^b_{\text{I}\R-c}(\text{I}(\K_X))$. We refer \cite{KS1} and
\cite{L} for a detailed study.

\section{Functorial properties of the microsupport and regularity for ind-sheaves}

\hspace*{\parindent}The microsupport $SS(F)$ of an ind-sheaf $F$
was defined in \cite{KS2} by several equivalent definitions which
we don't recall here (see Lemma 4.1 of \cite{KS2}). We start this
section with a new equivalent definition, which will be usefull in
the study of the functorial properties of this geometric object.

Let $\gamma$ be a closed convex proper cone in $\R^n$ and
$X\subset \R^n$ be an open subset. Denote by $q_1, q_2:X\times
X\to X$ the first and second projections and denote by $s:X\times
X\to X$ the map $(x,y)\mapsto x-y$. For an open subset $W\subset
X$, we define the functor $\Phi_{\gamma, W}:
D^b(\text{I}(\K_{X}))\to D^b(\text{I}(\K_{X}))$ by setting: $$
\Phi_{\gamma, W}(F)=Rq_{1 !!}(\K_{s^{-1}(\gamma)\cap
q_1^{-1}(W)\cap q_2^{-1}(W)}\otimes q_2^{-1}F).$$ We shall write
$\Phi_\gamma$ instead of $\Phi_{\gamma, X}$.

Let us denote by $j$ the embedding $s^{-1}(\gamma)\to X\times X$
and set $\widetilde{q}_i=q_i\circ j$, for $i=1,2$. Note that
$\Phi_\gamma(F)\simeq Rq_{1 !!}(Rj_{!!}j^{-1}q_2^{-1}F)\simeq
R\widetilde{q_1}_{!!}\widetilde{q_2}^{-1}F$, for all $F\in
D^b(\text{I}(\K_{X}))$.

\begin{lemma}\label{SP:2}
Let $X$ be an open subset of $\R^n$, $p=(x_0;\xi_0)\in T^*X$ and
$F\in D^b(\mathrm{I}(\K_X))$. Then $p\notin SS(F)$ if and only if
there exist $F'\in D^b(\mathrm{I}(\K_X))$, with compact support
and isomorphic to $F$ in an open neighborhood of $x_0$, and a
closed convex proper cone $\gamma$ in $\R^n$, with $\xi_0\in
\mathrm{Int}(\gamma^{\circ})$, such that $R\widetilde{q}_{2 *
}\widetilde{q}_1^!F'\simeq 0$, in a neighborhood of $x_0$.
\end{lemma}

\begin{proof}[\textbf{Proof.}]
Suppose $p\notin SS(F)$. By Lemma 4.1 of \cite{KS2}, there exists
a conic open neighborhood $U$ of $p$ in $T^*X$ such that, for any
$G\in D^b(\K_X)$ with $\mathrm{supp}(G)\subset\subset \pi(U)$ and
$SS(G)\subset U\cup T_X^*X$, one has
$\mathrm{Hom}_{D^b(\mathrm{I}(\K_{X}))}(G,F)=0$. Let $W$ be a
relatively compact open neighborhood of $x_0$, $\gamma$ a closed
convex proper cone such that $\xi_0\in\text{Int}(\gamma^{\circ
})$, $W\times \gamma^{\circ }\subset U\cup T_X^*X$ and choose
$\gamma^a$-open sets $\Omega_0\subset\Omega_1$, with
$\Omega_1\backslash\Omega_0\subset\subset W$ and
$x_0\in\text{Int}(\Omega_1\backslash\Omega_0)$. Let $I$ be a small
and filtrant category and $I\to D^b(\K_X);i\mapsto F_i$ a functor
such that $J(F)\simeq \underset{\underset{i\in I
}{\longrightarrow}}{\sideset{``}{"}\lim }J(F_i)$.

For each $j\in I$, set $H_j=R\widetilde{q}_{2 *
}\widetilde{q}_1^!R\Gamma_{\Omega_1\backslash\Omega_0}(F_j)\otimes\K_W$
and $G_j=R\widetilde{q}_{1 !}
\widetilde{q}_2^{-1}H_j\otimes\K_{\Omega_1\backslash\Omega_0}$. By
Corollary 6.4.5 and Proposition 5.2.3 of \cite{KS3}, one has
$SS(G_j)\subset\overline{\Omega_1\backslash\Omega_0}\times
\gamma^{\circ }$, for all $j\in I$. Then, for each $j\in I$, we
get, by assumption and the adjunction formulas:
$$0=\text{Hom}_{D^b(\mathrm{I}(\K_{X}))}(G_j,F)\simeq
\underset{\underset{i\in I }{\longrightarrow}}{\lim
}\text{Hom}_{D^b(\K_{X})}(G_j,F_i)\simeq$$ $$\simeq
\underset{\underset{i\in I }{\longrightarrow}}{\lim
}\text{Hom}_{D^b(\K_{X})}(H_j,R\widetilde{q}_{2 *
}\widetilde{q}_1^!R\Gamma_{\Omega_1\backslash\Omega_0}(F_i)).$$ It
follows that, for each $j\in I$, there exists $j\to i$ such that:
$$R\widetilde{q}_{2 *
}\widetilde{q}_1^!R\Gamma_{\Omega_1\backslash\Omega_0}(F_j)|_W\to
Rq_{2
* }\widetilde{q}_1^!R\Gamma_{\Omega_1\backslash\Omega_0}(F_i)|_W,$$
is the zero morphism. Since, for all $k\in\Z$:
$$H^k(R\widetilde{q}_{2
* }\widetilde{q}_1^!R\mathcal{I}{hom}(\K_{\Omega_1\backslash\Omega_0},F))\simeq
\underset{\underset{i\in I }{\longrightarrow}}{\sideset{``}{"}\lim
}H^k(R\widetilde{q}_{2 *
}\widetilde{q}_1^!R\Gamma_{\Omega_1\backslash\Omega_0}(F_i)),$$ we
get the desired result by taking
$F'=R\mathcal{I}{hom}(\K_{\Omega_1\backslash\Omega_0},F)$.

Coversely, let $W$ be an open neighborhood of $x_0$ such that
$F|_W\simeq F'|_W$ and $(R\widetilde{q}_{2
*}\widetilde{q}_1^!F')|_W\simeq 0$. For each $G\in D^b(\K_X)$,
with $\mathrm{supp}(G)\subset\subset W$ and $SS(G)\subset
W\times\gamma^{\circ}$, one has $G\simeq R\widetilde{q}_{1!}
\widetilde{q}_2^{-1}G_W$, by Proposition 5.2.3 of \cite{KS3}.
Hence: $$\text{Hom}_{D^b(\mathrm{I}(\K_{X}))}(G,F')\simeq
\text{Hom}_{D^b(\mathrm{I}(\K_{X}))}(R\widetilde{q}_{1!}
\widetilde{q}_2^{-1}G_W,F')\simeq$$
$$\simeq\text{Hom}_{D^b(\mathrm{I}(\K_{X}))}(G_W,R\widetilde{q}_{2
*}\widetilde{q}_1^!F')\simeq 0.$$ This implies that $p\notin
SS(F)$, by Lemma 4.1 of \cite{KS2}.
\end{proof}

Let $X$ and $Y$ be two real manifolds. Let us now denote by $q_1$
and $q_2$ the first and second projection defined on $X\times Y$.

Let $F\in D^b(\text{I}(\K_X))$ and $G\in D^b(\text{I}(\K_Y))$. One
sets $$F\boxtimes G=q_1^{-1}F\otimes q_2^{-1}G,$$ and call this
ind-sheaf the \emph{external tensor product} of $F$ and $G$.

\begin{proposition}\label{SP:4}
Let $F\in D^b(\mathrm{I}(\K_X))$ and $G\in D^b(\mathrm{I}(\K_Y))$.
Then: $$SS(F\boxtimes G) = SS(F)\times SS(G).$$
\end{proposition}

\begin{proof}[\textbf{Proof.}]
We may assume $X$ and $Y$ are vector spaces. Let $(p,p')\notin
SS(F)\times SS(G)$, with $p=(x_0;\xi_0)$ and $p'=(y_0;\eta_0)$.
Assume, for example, that $p\notin SS(F)$. Let $\gamma'$,
$\varepsilon>0$ and $W$ satisfying the condition $(4a)$ in Lemma
4.1 of \cite{KS2} for the ind-sheaf $F$. Set $\gamma=\gamma'\times
\{0\}$ and let $W'$ be a relatively compact open neighborhood of
$y_0$.

Let $s:(X\times Y)\times (X\times Y)\to X\times Y$ be the map
$(x,y,x',y')\mapsto(x-x',y-y')$, $j$ the embedding
$s^{-1}(\gamma)\hookrightarrow  X\times Y$ and let
$\widetilde{q}_1,\widetilde{q}_2:s^{-1}(\gamma) \to X\times Y$
denote the first and second projections. One has:
$$\Phi_{\gamma}(F_W\boxtimes G_{W'} )\simeq
R\widetilde{q_1}_{!!}\widetilde{q_2}^{-1}(F_W\boxtimes
G_{W'})\simeq$$ $$\simeq
R\widetilde{q_1}_{!!}(\widetilde{q_2}^{-1}q_1^{-1}(F_W)\otimes
\widetilde{q_2}^{-1}q_2^{-1}(G_{W'}))\simeq $$ $$\simeq
R\widetilde{q_1}_{!!}(\widetilde{q_2}^{-1}q_1^{-1}(F_W)\otimes
\widetilde{q_1}^{-1}q_2^{-1}(G_{W'}))\simeq$$ $$\simeq
R\widetilde{q_1}_{!!}(\widetilde{q_2}^{-1}q_1^{-1}(F_W))\otimes
q_2^{-1}(G_{W'}),$$ since $\widetilde{q_2}^{-1}$ commutes with
$\otimes$, $q_2\circ\widetilde{q_2}=q_2\circ\widetilde{q_1}$ and
applying Proposition 5.2.7 of \cite{KS1}. Let now $t$ denote the
map $X\times X\to X;(x,x')\mapsto x-x'$, $\widetilde{p_1},
\widetilde{p_2}:t^{-1}(\gamma') \to X$ denote the first and second
projections and let $k:s^{-1}(\gamma)\to t^{-1}(\gamma')$ be the
map $(x,y,x',y')\mapsto(x,x')$. Then:
$$R\widetilde{q_1}_{!!}(\widetilde{q_2}^{-1}q_1^{-1}(F_W))\simeq
R\widetilde{q_1}_{!!}k^{-1}\widetilde{p_2}^{-1}(F_W)\simeq
q_1^{-1}R\widetilde{p_1}_{!!}\widetilde{p_2}^{-1}(F_W),$$ where
the first isomorphism follows from the equality $q_1\circ
\widetilde{q}_2=\widetilde{p}_2\circ k$ and the second from
Proposition 5.2.9 of \cite{KS1}. Since
$R\widetilde{p_1}_{!!}\widetilde{p_2}^{-1}(F_W)\simeq
\Phi_{\gamma'}(F_W)$, we obtain: $$\Phi_{\gamma, W\times
W'}(F\boxtimes G)\simeq(\Phi_{\gamma}(F_W\boxtimes
G_{W'}))_{W\times W'}\simeq$$ $$\simeq
(\Phi_{\gamma'}(F_W))_W\boxtimes G_{W'}\simeq \Phi_{\gamma',
W}(F)\boxtimes G_{W'}\simeq 0,$$ which implies
$(x_0,y_0;\xi_0,\eta_0)\notin SS(F\boxtimes G)$, by Lemma 4.1 of
\cite{KS2}.

Conversely, let $(p,p')\in SS(F)\times SS(G)$. By Lemma 4.1 of
\cite{KS2}, we may find a small and filtrant category $I$,
integers $a\leq b$ and a functor $I\to
C^{[a,b]}(\mathrm{Mod}(\K_X));i\mapsto F_i$, such that $J(F)\simeq
\underset{\underset{i\in I }{\longrightarrow}}{\sideset{``}{"}\lim
}J(F_i)$ and, for all conic open neighborhood $U$ of $p$ in
$T^*X$, there exists $i\in I$ such that every morphism $i\to i'$
in $I$ induces a non-zero morphism $F_i\to F_{i'}$ in
$D^b(\K_X;U)$. Similarly, we may find a small and filtrant
category $J$, integers $a\leq b$ and a functor $J\to
C^{[a,b]}(\mathrm{Mod}(\K_Y));j\mapsto G_j$, such that $J(G)\simeq
\underset{\underset{j\in J }{\longrightarrow}}{\sideset{``}{"}\lim
}J(G_j)$ and, for all conic open neighborhood $V$ of $p'$ in
$T^*Y$, there exists $j\in J$ such that every morphism $j\to j'$
in $J$ induces a non-zero morphism $G_j\to G_{j'}$ in
$D^b(\K_Y;V)$.

Then $J(F\boxtimes G)\simeq \underset{\underset{i\in I, j\in J
}{\longrightarrow}}{\sideset{``}{"}\lim }J(F_i\boxtimes G_j)$ and
every morphism $(i,j)\to (i',j')$ in $I\times J$\footnote{We
denote by $I\times J$ the product of the categories $I$ and $J$.
See \cite{KS4}.} induces a non-zero morphism $F_i\boxtimes G_j\to
F_{i'}\boxtimes G_{j'}$ in $D^b(\K_{X\times Y};U\times V)$. Since
$U\times V$ forms a neighborhood system of $(p,p')$, we may
conclude that $(p,p')\in SS(F\boxtimes G)$.
\end{proof}

Let $f$ be a morphism from $Y$ to $X$. We denote by $f_d$ and
$f_\pi$ the canonical morphisms ( $f_d$ was noted by ${}^{t} f'$
in \cite{KS3}):
\begin{center}
$f_d:Y\times_X T^*X\rightarrow T^*Y$ and $f_\pi:Y\times_X
T^*X\rightarrow T^*X$.
\end{center}

\begin{proposition}\label{SP:5}
Let $M$ be a closed submanifold of $X$ and let $j$ denote the
embedding $M\hookrightarrow X$. Let $G\in
D^{b}(\mathrm{I}(\K_M))$. Then, $$SS(Rj_*G) =
j_{\pi}j_d^{-1}(SS(G)).$$
\end{proposition}

\begin{proof}[\textbf{Proof}]
Let $d$ denote the codimension of $M$. We may assume $X=\R^n$ and
$M=\{0\}\times \R^{n-d}$. Let $(x_1,...,x_n)$ be a system of
coordinates on $X$, set $x'=(x_1,...,x_d)$,
$x''=(x_{d+1},...,x_n)$ and let $(x',x'';\xi',\xi'')$ be the
associated coordinates on $T^*X$.

Let $(x'_0,x''_0;\xi'_0,\xi''_0)\notin j_{\pi}j_d^{-1}(SS(G))$. We
may assume $x'_0=0$. Hence, $(x''_0;\xi''_0)\notin SS(G)$ and we
may find a closed convex proper cone $\gamma$ in $\R^{n-d}$, with
$\xi''_0\in \text{Int}(\gamma^{\circ a})$, and an open
neighborhood $W$ of $(0,x''_0)$ in $M$ such that, for all $F\in
D^b(\K_M)$, with $\mathrm{supp}(F)\subset\subset W$ and
$SS(F)\subset W\times \gamma^{\circ a}$, one has
$\mathrm{Hom}_{D^b(\text{I}(\K_M))}(F,G)=0$. Let $V$ be an open
neighborhood of $(0,x''_0)$ in $X$ such that $W=V\cap M$. For each
$F\in D^b(\K_X)$ such that $\mathrm{supp}(F)\subset\subset V$ and
$SS(F)\subset V\times (\R^d\times\gamma^{\circ a})$, one has, by
Theorem 5.2.4 of \cite{KS1}:
$$\mathrm{Hom}_{D^b(\text{I}(\K_X))}(F,Rj_*G)\simeq\mathrm{Hom}_{D^b(\text{I}(\K_M))}(j^{-1}F,G).$$
Since $\mathrm{supp}(j^{-1}F)\subset\subset W$ and
$SS(j^{-1}F)\subset W\times \gamma^{\circ a}$, it follows from the
hypothesis that
$$\mathrm{Hom}_{D^b(\text{I}(\K_X))}(F,Rj_*G)\simeq 0.$$

Conversely, let $p=(x_0;\xi_0)\notin SS(Rj_*G)$. We may assume
$x_0\in M$. Take $\gamma$ and $W$ satisfying the condition $(4a)$
for the ind-sheaf $Rj_*G$. By Proposition 5.2.9 of \cite{KS1}, we
may prove that $$0\simeq j^{-1}\Phi_{\gamma, W}(Rj_*G)\simeq
\Phi_{\gamma\cap M, W\cap M}(G).$$ From this we get
$(x''_0;\xi''_0)\notin SS(G)$.
\end{proof}

\begin{proposition}\label{SP:6}
Let $Y$, $X$ be real analytic manifolds, $f\colon Y\to X$ be a
smooth morphism and $F\in D^b (\mathrm{I}(\K_{X}))$. Then
$$SS(f^{-1}F) = f_d f_{\pi}^{-1}(SS(F).$$
\end{proposition}

\begin{proof}[\textbf{Proof}]
We may assume $Y=X\times Z$, $f$ is the first projection and that
$Y, X$ and $Z$ are vector spaces. Then, $f^{-1}F\simeq F\boxtimes
\K_Z$ and the desired equality follows from Proposition
\ref{SP:4}.
\end{proof}

\begin{proposition}\label{SP:7}
Let $F\in D^b(\mathrm{I}(\K_X))$ and $G\in D^b(\mathrm{I}(\K_Y))$.
Then: $$SS(R\mathcal{I}{hom}(q_2^{-1}G,q_1^{!}F))\subset
SS(F)\times SS(G)^a.$$
\end{proposition}

\begin{proof}[\textbf{Proof}]
We may assume $X$ and $Y$ are vector spaces. Let
$p=(x_0,y_0;\xi_0,\eta_0)\notin SS(F)\times SS(G)^a$.

First assume $(x_0;\xi_0)\notin SS(F)$. By Lemma \ref{SP:2}, there
exist a proper closed convex cone $\gamma'$, with
$\xi\in\text{Int}(\gamma'^\circ)$, and $F'\in
D^b(\mathrm{I}(\K_X))$ with compact support and isomorphic to $F$
in a neighborhood of $x_0$ such that
$R\widetilde{p_2}_*\widetilde{p_1}^!F'\simeq 0$ in a neighborhood
of $x_0$, where $t$ denotes the map $X\times X\to X;(x,x')\mapsto
x-x'$ and $\widetilde{p_1}, \widetilde{p_2}:t^{-1}(\gamma') \to X$
denote the first and second projections.

Set $\gamma:=\gamma'\times\{0\}$ and let $W'$ be a relatively
compact open neighborhood of $y_0$. Denote by $s:(X\times Y)\times
(X\times Y)\to X\times Y$ the map $(x,y,x',y')\mapsto
(x-x',y-y')$, by $j$ the embedding $s^{-1}(\gamma)\hookrightarrow
(X\times Y)\times(X\times Y)$, by $\widetilde{q_1},
\widetilde{q_2}:s^{-1}(\gamma)\to X\times Y$ the first and second
projections and let $k$ denote the map $(X\times Y)\times(X\times
Y)\to X\times X; (x,y,x',y')\mapsto (x,x')$. We shall prove that:
\begin{equation}\label{E:E}
R\widetilde{q_2}_{*}\widetilde{q_1}^!R\mathcal{I}{hom}(q_2^{-1}G_{W'},q_1^{!}F')\simeq
0,\ \text{in a neighborhood of $(x_0,y_0)$.}
\end{equation}

By Propositions 5.3.8 and 5.2.3 of \cite{KS1}, one has:
$$R\widetilde{q_2}_{*}\widetilde{q_1}^!R\mathcal{I}{hom}(q_2^{-1}G_{W'},q_1^{!}F')\simeq
R\widetilde{q_2}_{*}R\mathcal{I}{hom}(\widetilde{q_1}^{-1}q_2^{-1}G_{W'},\widetilde{q_1}^!q_1^{!}F')\simeq$$
$$\simeq
R\widetilde{q_2}_{*}R\mathcal{I}{hom}(\widetilde{q_2}^{-1}q_2^{-1}G_{W'},\widetilde{q_1}^!q_1^{!}F')\simeq
R\mathcal{I}{hom}(q_2^{-1}G_{W'},R\widetilde{q_2}_{*}\widetilde{q_1}^!q_1^{!}F'),$$
and $$R\widetilde{q_2}_{*}\widetilde{q_1}^!q_1^{!}F'\simeq
R\widetilde{q_2}_{*}(k\circ j)^!\widetilde{p_1}^!F'\simeq
q_1^!R\widetilde{p_2}_*\widetilde{p_1}^!F',$$ where the last
isomorphism follows from Proposition 5.3.10 of \cite{KS1}. By
assumption, $q_1^!R\widetilde{p_2}_*\widetilde{p_1}^!F'\simeq 0$,
in a neighborhood of $(x_0,y_0)$, which entails (\ref{E:E}).

Now assume $(y_0;\eta_0)\notin SS(G)$. Take a cone $\gamma'$, with
$\eta_0\in \text{Int}(\gamma^\circ)$, and $G'\in D^b(\K_{X_{sa}})$
with compact support and isomorphic to $G$ in a neighborhood of
$y_0$ such that $\Phi_{\gamma'}(G')\simeq 0$ in a neighborhood of
$y_0$. Set $\gamma:=\{0\}\times\gamma'^a$ and let $W$ be a
relatively compact open neighborhood of $x_0$. Denote by
$s:(X\times Y)\times (X\times Y)\to X\times Y$ the map
$(x,y,x',y')\mapsto (x-x',y-y')$, by $j$ the embedding
$s^{-1}(\gamma)\hookrightarrow (X\times Y)\times(X\times Y)$ and
by $\widetilde{q_1}, \widetilde{q_2}:s^{-1}(\gamma)\to X\times Y$
the first and second projections. Let us prove that
$$R\widetilde{q_2}_{*}\widetilde{q_1}^!R\mathcal{I}{hom}(q_2^{-1}G',q_1^{!}F_W)\simeq
0,$$ in a neighborhood of $(x_0,y_0)$.

Let now $k:(X\times Y)\times (X\times Y)\to Y\times Y$ be the map
$(x,y,x',y')\mapsto(y',y)$, $t$ denote the map $Y\times Y\to
Y;(y,y')\mapsto y-y'$ and $\widetilde{p_1},
\widetilde{p_2}:t^{-1}(\gamma') \to Y$ denote the first and second
projections. By Propositions 5.3.8 and 5.3.5 of \cite{KS1}, one
has:
$$R\widetilde{q_2}_{*}\widetilde{q_1}^!R\mathcal{I}{hom}(q_2^{-1}G',q_1^{!}F_W)\simeq
R\widetilde{q_2}_{*}R\mathcal{I}{hom}(\widetilde{q_1}^{-1}q_2^{-1}G',\widetilde{q_1}^!q_1^{!}F_W)\simeq$$
$$\simeq
R\widetilde{q_2}_{*}R\mathcal{I}{hom}(\widetilde{q_1}^{-1}q_2^{-1}G',\widetilde{q_2}^!q_1^{!}F_W)\simeq
R\mathcal{I}{hom}(R\widetilde{q_2}_{!!}\widetilde{q_1}^{-1}q_2^{-1}G',q_1^{!}F_W),$$
and by Proposition 5.2.9 of \cite{KS1}:
$$R\widetilde{q_2}_{!!}\widetilde{q_1}^{-1}q_2^{-1}G'\simeq
R\widetilde{q_2}_{!!}(k\circ j)^{-1}\widetilde{p_2}^{-1}G'\simeq
q_2^{-1}R\widetilde{p_1}_{!!}\widetilde{p_2}^{-1}G'\simeq
q_2^{-1}\Phi_{\gamma'}(G').$$ The result follows.
\end{proof}

Recall that is possible to localize the category
$D^b(\mathrm{I}(\K_X))$ with respect to the microsupport as for
classical sheaves. In fact, given $V\subset T^*X$, one sets:
$D^b(\mathrm{I}(\K_X;\Omega))=D^b(\mathrm{I}(\K_X))/D_V^b(\mathrm{I}(\K_X)),$
where $\Omega=T^*X\backslash V$ and $D_V^b(\mathrm{I}(\K_X))$ is
the full triangulated subcategory of $D^b(\mathrm{I}(\K_X))$
consisting of objects $F$ such that $SS(F)\subset V$.

We shall now prove the counterpart for ind-sheaves of the
microlocal characterization of classical sheaves whose
microsupport is contained in an involutive manifold.

Let $f:Y\to X$ be a morphism of manifolds.

\begin{proposition}\label{SP:8}
Assume $f$ is a closed embedding and identify $Y$ with a
submanifold of $X$.

(i) Let $G\in D^b(\mathrm{I}(\K_Y))$. Then $Rf_*G$ is regular
along $\pi^{-1}(Y)$.

(ii) Let $F\in D^b(\mathrm{I}(\K_X))$ and assume $F$ is regular
along $\pi^{-1}(Y)$ at $p\in T_Y^*X$. Then there exists $G\in
D^b(\mathrm{I}(\K_Y))$ such that $F\simeq Rf_*G$ in
$D^b(\mathrm{I}(\K_X;p))$.
\end{proposition}

\begin{proof}[\textbf{Proof}]
(i) Consider a small and filtrant category $I$ and a functor $I\to
C^{[a,b]}(\K_Y);i\mapsto G_i$ such that $G\simeq
Q(\underset{\underset{i\in I
}{\longrightarrow}}{\sideset{``}{"}\lim }G_i)$. Since $f$ is
proper on the support of $G$ one has, by Theorem 3.1 of \cite{KS2}
together with Proposition 2.3.2 of \cite{KSFI}, $J(Rf_*G)\simeq
\underset{\underset{i\in I }{\longrightarrow}}{\sideset{``}{"}\lim
}J(Rf_*G_i)$ and, for each $i\in I$, $$SS(Rf_*G_i)\subset
\pi^{-1}(Y).$$ It follows that $Rf_*G$ is regular along
$\pi^{-1}(Y)$ at each $p\in T^*X$.

(ii) The proof is an adaptation of the proof of Proposition 6.6.1
of \cite{KS3}.

By hypothesis, there exist $F'$ isomorphic to $F$ in a
neighborhood of $\pi(p)$, a conic open neighborhood $U$ of $p$, a
small and filtrant category $I$ and a functor $I\to
D^{[a,b]}(\K_X);i\mapsto F_i$ such that $J(F')\simeq
\underset{\underset{i\in I }{\longrightarrow}}{\sideset{``}{"}\lim
}J(F_i)$ and $SS(F_i)\cap U\subset \pi^{-1}(Y)$, for all $i\in I$.
We may assume from the beginning that $J(F)\simeq
\underset{\underset{i\in I }{\longrightarrow}}{\sideset{``}{"}\lim
}J(F_i)$.

If $p\in T_X^*X$, then there exists an open neighborhood $W$ of
$\pi(p)$ such that $\pi^{-1}(W)\subset U$ and
\begin{equation}\label{SE:1}
F_i|_W\simeq (F_i)_Y|_W, \ \text{for all $i\in I$.}
\end{equation}

One has the following distinguished triangle in
$D^b(\text{I}(\K_X))$: $$F\otimes\K_{X\backslash Y}\to F\to
F\otimes\K_Y\xrightarrow{+1},$$ with $J(F\otimes\K_{X\backslash
Y})\simeq \underset{\underset{i\in I
}{\longrightarrow}}{\sideset{``}{"}\lim }J((F_i)_{X\backslash
Y})$, by Proposition 2.3.2 of \cite{KSFI}, and $SS((F_i)_{
X\backslash Y})\cap \pi^{-1}(W)=\emptyset$, by (\ref{SE:1}), for
each $i\in I$. Hence, $SS(F\otimes \K_{X\backslash Y})\cap
\pi^{-1}(W)=\emptyset$, by Lemma 4.1 of \cite{KS2}. It follows
that $F\simeq F\otimes\K_Y\simeq Rj_*j^{-1}F$ in
$D^b(\text{I}(\K_X;p))$.

Now assume $p\in \dot{T}_Y^*X$. We shall argue by induction on the
codimension of $Y$ and prove that there exist a conic open
neighborhood $U'\subset U$ of $p$, $G\in D^b(\text{I}(\K_X))$ such
that $F\simeq G\otimes\K_Y$ in $D^b(\text{I}(\K_X;U'))$ and
$J(G)\simeq\underset{\underset{i\in I
}{\longrightarrow}}{\sideset{``}{"}\lim }J(G_i)$, for some $G_i\in
D^b(\K_X)$ such that $SS(G_i\otimes\K_Y)\cap U'=SS(F_i)\cap U'$,
for all $i\in I$.

Let us first assume $Y$ is a hypersurface. Let $(x_1,...,x_n)$ be
a system of local coordinates of $X$ in a neighborhood of $\pi(p)$
such that $$Y=\{(x_1,...,x_n); x_1=0\},$$ and let
$(x_1,...,x_n;\xi_1,...,\xi_n)$ denote the associated coordinates
in $T^*X$. Set $x'=(x_2,...,x_n)$, $x=(x_1,x')$,
$\xi'=(\xi_2,...,\xi_n)$, $\xi=(\xi_1,\xi')$,
$\Omega^\pm=\{(x_1,...,x_n);$ $\pm x_1>0\}$ and denote by $j_\pm$
the open embeddings $\Omega^\pm\hookrightarrow X$. We may assume
$p=(0;1,0,...,0)$ and that $U=W\times\text{Int}\gamma$, where $W$
is an open neighborhood of $0$ and $\gamma$ is the cone
$\{(\xi_1,\xi'); \xi_1\geq \delta|\xi'|\}$, for some $\delta>0$.

Let us consider the following distinguished triangle in
$D^b(\text{I}(\K_X))$:
\begin{equation}\label{SE:2}
R\mathcal{I}{hom}(\K_{\{x_1\geq 0\}},F)\to F\to
R\mathcal{I}{hom}(\K_{\{x_1< 0\}},F)\xrightarrow{+1}.
\end{equation}
By Proposition 3.8 of \cite{KS2}, we have
$$J(R\mathcal{I}{hom}(\K_{\{x_1<
0\}},F))\simeq\underset{\underset{i}{\longrightarrow}}{\sideset{``}{"}\lim}J(R\mathcal{I}{hom}(\K_{\{x_1<
0\}},F_i))\simeq$$
$$\simeq\underset{\underset{i}{\longrightarrow}}{\sideset{``}{"}\lim}J(R\mathcal{H}{om}(\K_{\{x_1<
0\}},F_i))\simeq
\underset{\underset{i}{\longrightarrow}}{\sideset{``}{"}\lim}J(Rj_{-
*}j_-^{-1}(F_i)).$$

Let us prove that $U\cap SS(Rj_{- *}j_-^{-1}F_i)=\emptyset,$ for
all $i\in I$.

By Theorem 6.3.1 of \cite{KS3}, $$SS(Rj_{- *}j_-^{-1}F_i)\subset
SS(j_-^{-1}F_i)\widehat{+}N^*(\Omega^-),\ \text{for all $i\in
I$}.$$

Let $i\in I$ and assume that there exists $(x_0;\xi_0)\in U\cap
(SS(j_-^{-1}F_i)\widehat{+}N^*(\Omega^-))$. Let
$\{(x_{n};\xi_n)\}_n$ (resp. $\{(y_n;\eta_n)\}_n$) be a sequence
in $SS(j_-^{-1}F_i)$ (resp. in $N^*(\Omega^-)$), such that:
$$\begin{cases} x_n,y_n \xrightarrow[n]{}x_0,\\ \xi_n+\eta_n
\xrightarrow[n]{} \xi_0,\\ |x_n-y_n||\xi_n|\xrightarrow[n]{} 0.
\end{cases}$$

Since $U\cap SS(j_-^{-1}F_i)\subset \pi^{-1}(Y) \cap
T^*\Omega^-=\emptyset $, one has $(x_n;\xi_n)\notin U$, for each
$n\in\N$. On the other hand, since the sequence $\{x_n\}_n$
converges to $x_0$, $x_n\in \pi(U)$, for $n\gg 1$, and hence,
$\xi_n\notin \text{Int}\gamma$, for $n\gg 1$. Moreover, if the
sequence $\{\eta_n\}_n$ converge to $0$, then the sequence
$\{(x_n;\xi_n)\}_n$ will converge to $(x_0;\xi_0)$ and we must
have $(x_n;\xi_n)\in U$, for a sufficiently large $n$, which is a
contradiction. Hence, after replacing the sequences by convenient
subsequences, we may assume $\xi_n\notin \text{Int}\gamma$,
$y_n\in
\partial \Omega^-$, $\eta_n'=0$ and $\eta_{n, 1}\leq 0$, for all
$n\in\N$. From this we get that:
$$\xi_{0,1}=\underset{\underset{n}{\longrightarrow}}{\lim}(\xi_{n,1}+\eta_{n,1})
\leq\underset{\underset{n}{\longrightarrow}}{\lim}\delta|\xi'_n|=\delta|\xi'_0|,$$
which is a contradiction. It follows that $SS(R\Gamma_{\{x_1\geq
0\}}(F_i))\cap U=SS(F_i)\cap U$, for all $i\in I$, and
$SS(R\mathcal{I}{hom}(\K_{\{x_1< 0\}},F))\cap U=\emptyset$.
Moreover, $F\simeq R\mathcal{I}{hom}(\K_{\{x_1\geq 0\}},F)$ in
$D^b(\text{I}(\K_X, U))$.

Set $G=R\mathcal{I}{hom}(\K_{\{x_1\geq 0\}},F)$. One has:
$$J(G\otimes\K_{X\backslash
Y})\simeq\underset{\underset{i}{\longrightarrow}}{\sideset{``}{"}\lim}J(Rj_{+
!}j_+^{-1}(R\Gamma_{\{x_1\geq 0\}}(F_i))),$$ and we may prove as
above that $U\cap SS(Rj_{+ !}j_+^{-1}(R\Gamma_{\{x_1\geq
0\}}(F_i)))=\emptyset$, for all $i\in I$, which entails
$SS(G\otimes\K_{X\backslash Y})\cap U=\emptyset$ and
$SS(R\Gamma_{\{x_1\geq 0\}}(F_i)\otimes\K_Y)\cap U=SS(F_i)\cap U$,
for all $i\in I$. Moreover, $F\simeq G\otimes \K_Y$ in
$D^b(\text{I}(\K_X, U))$.

Let us now assume that $Y$ is a submanifold of $X$ of codimension
$m\geq 1$ and let $(x_1,...,x_n)$ be a system of local coordinates
of $X$ in a neighborhood of $\pi(p)$ such that $Y=\{(x_1,...,x_n);
x_1=...=x_m=0\}$ and let $(x_1,...,x_n;\xi_1,...,\xi_n)$ denote
the associated coordinates in $T^*X$. We may assume
$p=(0;1,0,...0)$. Let $Y'$ be the submanifold of $X$:
$$\{(x_1,...,x_n)\in X; x_1=...=x_{m-1}=0\},$$ and let us denote
by $\iota$ (resp. $k$) the embedding $Y'\hookrightarrow X$ (resp.
$Y\hookrightarrow Y'$).

Since $\pi^{-1}(Y)\subset\pi^{-1}(Y')$, by induction, there exist
a conic open neighborhood $U'\subset U$ of $p$, $G\in
D^b(\text{I}(\K_X))$ such that $G\otimes \K_{Y'}\simeq F$ in
$D^b(\text{I}(\K_X;U'))$ and $J(G)\simeq\underset{\underset{i\in I
}{\longrightarrow}}{\sideset{``}{"}\lim }J(G_i)$, for some $
G_i\in D^b(\K_X)$ with  $SS(G_i\otimes\K_{Y'})\cap U'=SS(F_i)\cap
U'$, for all $i\in I$.  Since
$$SS(G_i\otimes\K_{Y'})=SS(\iota_*\iota^{-1}G_i)=\iota_\pi\iota_d^{-1}(SS(\iota^{-1}G_i)),$$
we get, $SS(\iota^{-1}G_i)\cap\iota_d\iota_\pi^{-1}(U')\subset
\pi^{-1}(Y)$, for all $i\in I$, where $\iota_d\iota_\pi^{-1}(U')$
is a conic open neighborhood of $p'=(0;1,0,...,0)\in
\dot{T}_Y^*Y'$.

Since $Y$ is a hypersurface in $Y'$, we may find a conic open
neighborhood $V\subset \iota_d\iota_\pi^{-1}(U')$ of $p'$, $H\in
D^b(\text{I}(\K_{Y'}))$ in $D^b(\text{I}(\K_{Y'};V))$ such that
$H\otimes \K_Y\simeq \iota^{-1}G$ and
$J(H)\simeq\underset{\underset{i\in I
}{\longrightarrow}}{\sideset{``}{"}\lim }J(H_i)$, for some $H_i\in
D^b(\K_{Y'})$ with $SS(H_i\otimes\K_Y)\cap V=SS(\iota^{-1}G_i)\cap
V$, for all $i\in I$. We may assume $V=(V_1\cap Y')\times V_2$,
for some open neighborhood $V_1$ (resp. $V_2$) of $0$ in $X$
(resp. of $(1,...,0)$ in $\R^{n-m+1}$). Set $V'=V_1\times
(\R^{m-1}\times V_2)\subset T^*X$. It follows by Proposition
\ref{SP:5} that $\iota_*H\otimes\K_Y\simeq G\otimes\K_{Y'}$.
Moreover, $SS(\iota_*H_i\otimes\K_Y)\cap
V'=SS(\iota_*(H_i\otimes\K_Y))\cap
V'=\iota_\pi\iota_d^{-1}(SS(H_i\otimes\K_Y)\cap
V)=\iota_\pi\iota_d^{-1}(SS(\iota^{-1}G_i)\cap
V)=SS(\iota_*\iota^{-1}G_i)\cap V'=SS(G_i\otimes\K_{Y'})\cap V'$.
The result follows.
\end{proof}

\begin{proposition}\label{SP:9}
Assume $f$ is smooth and identify $Y\times_ X T^*X$ with a
submanifold of $T^*Y$.

(i) Let $F\in D^b(\mathrm{I}(\K_X))$. Then $f^{-1}F$ is regular
along $Y\times_ X T^*X$.

(ii) Let $G\in D^b(\mathrm{I}(\K_Y))$, $p\in Y\times_ X T^*X$ and
assume $G$ is regular along $Y\times_ X T^*X$ at $p$. Then there
exists $F\in D^b(\mathrm{I}(\K_X))$ such that $G\simeq f^{-1}F$ in
$D^b(\mathrm{I}(\K_Y;p))$.
\end{proposition}

\begin{proof}[\textbf{Proof}]
(i) Consider a small and filtrant category $I$ and a functor $I\to
C^{[a,b]}(\K_X);i\mapsto F_i$ such that $F\simeq
Q(\underset{\underset{i\in I
}{\longrightarrow}}{\sideset{``}{"}\lim }F_i)$. Then, Theorem 3.2
of \cite{KS2} together with Proposition 2.3.2 \cite{KSFI} entails
$J(f^{-1}F)\simeq \underset{\underset{i\in I
}{\longrightarrow}}{\sideset{``}{"}\lim }J(f^{-1}F_i)$ and one has
$$SS(f^{-1}F_i)\subset Y\times_ X T^*X, \ \text{for all $i\in
I$}.$$ It follows that $f^{-1}F$ is regular along $Y\times_ X
T^*X$ at each $p\in T^*Y$.

(ii)  The proof is an adaptation of Proposition 6.6.2 of
\cite{KS3}.

By the assumption, there exist $G'$ isomorphic to $G$ in a
neighborhood of $\pi(p)$, a conic open neighborhood $U$ of $p$, a
small and filtrant category $I$ and a functor $I\to
D^{[a,b]}(\K_Y);i\mapsto G_i$ such that $J(G')\simeq
\underset{\underset{i\in I }{\longrightarrow}}{\sideset{``}{"}\lim
}J(G_i)$ and $SS(G_i)\cap U\subset Y\times_ X T^*X$, for all $i\in
I$. Since $J(G_i)\simeq\underset{\underset{U\subset\subset Y
}{\longrightarrow}}{\sideset{``}{"}\lim }J((G_i)_U)$, where $U$
runs through the family of relatively compact open subsets of $Y$,
we may assume from the beginning that $J(G)\simeq
\underset{\underset{i\in I }{\longrightarrow}}{\sideset{``}{"}\lim
}J(G_i)$ and $G_i$ has compact support, for all $i\in I$. We may
also assume that $Y=Z\times X$, that $f$ is the projection
$Z\times X\to X$ and that $Y, Z, X$ are vector spaces.

We shall prove, by induction on $\text{dim} Z$, that there exist a
conic open neighborhood $U'\subset U$ of $p$ and $F\in
D^b(\text{I}(\K_X))$ such that $f^{-1}F\simeq G$ in
$D^b(\text{I}(\K_Y; U'))$ and $J(F)\simeq\underset{\underset{i\in
I }{\longrightarrow}}{\sideset{``}{"}\lim }J(F_i)$, for some $F_i
\in D^b(\K_X)$, with $SS(f^{-1}F_i)\cap U'=SS(G_i)\cap U'$, for
all $i\in I$.

Let us first suppose that dim$Z=1$. We may assume $Y=\R^n$,
$X=\R^{n-1}$, $p=(0;\xi_0)$ and $\xi_0=(0,\xi'_0)$. Denote by
$q_1, q_2:Y\times Y\to Y$ the first and second projections and
denote by $s:Y\times Y\to Y$ the map $(x,y)\mapsto x-y$.

If $\xi_0=0$, there exists an open neighborhood $W$ of $0$ such
that $\pi^{-1}(W)\subset U$. Let $\gamma$ be the closed cone
$\R\times\{0\}\subset \R^n$ and $\Omega_0\subset\Omega_1$ be two
$\gamma$-open sets such that
$\Omega_1\backslash\Omega_0\subset\subset W$ and
$0\in\text{Int}(\Omega_1\backslash\Omega_0)$. Set
$G'=G\otimes\K_{\Omega_1\backslash\Omega_0}$ and
$G'_i=G_i\otimes\K_{\Omega_1\backslash\Omega_0}$, for all $i\in
I$. One has $J(G')\simeq\underset{\underset{i\in I
}{\longrightarrow}}{\sideset{``}{"}\lim }J(G'_i)$ and, by
assumption, $SS(G'_i)\subset Y\times_X T^*X$, for all $i\in I$. By
the microlocal cut-off Lemma (Proposition 5.2.3 of \cite{KS3}),
one has $G'_i\simeq \phi_{\gamma}^{-1}R\phi_{\gamma *}G'_i$, for
all $i\in I$.

Let us consider the following distinguished triangle in
$D^b(\text{I}(\K_Y))$: $$Rq_{1 !!}(\K_{s^{-1}(\gamma)}\otimes
q_2^{-1}G')\to G'\to Rq_{1 !!}(K\otimes q_2^{-1}G')\xrightarrow
{+1},$$ where $K$ denotes the complex $\K_{s^{-1}(\gamma)}\to
\K_{s^{-1}(0)}$, with $\K_{s^{-1}(0)}$ in degree $0$.

Since $q_1$ is proper on the support of $q_2^{-1}G'_i$, for all
$i\in I$, one has by Proposition 2.3.2 of \cite{KSFI}: $$J(Rq_{1
!!}(\K_{s^{-1}(\gamma)}\otimes q_2^{-1}G'))\simeq
\underset{\underset{i}{\longrightarrow}}{\sideset{``}{"}
\lim}J(Rq_{1 !}(\K_{s^{-1}(\gamma)}\otimes q_2^{-1}G'_{i })),$$
and $$J(Rq_{1 !!}(K\otimes
q_2^{-1}G'))\simeq\underset{\underset{i}{\longrightarrow}}{\sideset{``}{"}
\lim} J(Rq_{1 !}(K\otimes q_2^{-1}G'_{i})).$$ Moreover, $Rq_{1
!}(\K_{s^{-1}(\gamma)}\otimes q_2^{-1}G'_{i }
)\simeq\phi_{\gamma}^{-1}R\phi_{\gamma *}G'_i\simeq G'_i,$ for all
$i\in I$, which entails $Rq_{1!}(K\otimes q_2^{-1}G'_i)\simeq 0$.
On the other hand, by Theorem 5.2.9 of \cite{KS1}, one has:
$$Rq_{1 !!}(\K_{s^{-1}(\gamma)}\otimes q_2^{-1}G')\simeq
f^{-1}Rf_{!!}G'.$$ Setting
$W'=\text{Int}(\Omega_1\backslash\Omega_0)$, we conclude that
$G|_{W'}\simeq f^{-1}Rf_{!!}G'|_{W'}$ and
$SS(f^{-1}Rf_{!}G_i')\cap\pi^{-1}(W')=SS(G'_i)\cap\pi^{-1}(W')=SS(G_i)\cap\pi^{-1}(W')$,
for all $i\in I$.

Now assume $\xi_0\neq 0$. Let $\varepsilon>0$ and set
$H_\varepsilon=\{x\in\R^n;\langle x;\xi_0\rangle>-\varepsilon\}$.
Let $\gamma$ be a proper closed convex cone in $\R^n$, with
$\xi_0\in\text{Int}(\gamma^{\circ a})$, and let $W$ be an open
neighborhood of $0$ which is the intersection of $H_\varepsilon$
and a $\gamma$-open subset of $\R^n$, for some $0<\varepsilon\ll
1$, such that $W\times \gamma^{\circ a}\backslash\{0\}\subset U$.
Then, for all $i\in I$: $$SS(G_i)\cap (W\times \gamma^{\circ
a})\subset Y\times_X T^*X.$$

Set $G'_i=\phi_{\gamma}^{-1}R\phi_{\gamma *}(G_{i,
H_\varepsilon})$, for all $i\in I$. Arguing as M. Kashiwara and P.
Schapira in the proof of Proposition 6.6.2 of \cite{KS3}, we may
prove that $SS(G_i')\cap (W\times \gamma^{\circ a})\subset
Y\times_X T^*X$, for all $i\in I$, which is equivalent, by
Proposition 5.2.3 of \cite{KS3}, to $SS(G'_i)\cap
\pi^{-1}(W)\subset Y\times_X T^*X$.

Let us now set $G'=Rq_{1 !!}(\K_{s^{-1}(\gamma^a)}\otimes
q_2^{-1}(G\otimes\K_{H_\varepsilon}))$. Arguing as in the previous
case $(\xi=0)$ we may prove that there exist an open neighborhood
$W'\subset W\cap H_\varepsilon$ of $0$ and $F\in
D^b(\text{I}(\K_X))$ such that $G'|_{W'}\simeq f^{-1}F|_{W'}$ and
$J(F)\simeq\underset{\underset{i\in I
}{\longrightarrow}}{\sideset{``}{"}\lim }J(F_i)$, for some $F_i\in
D^b(\K_X)$ with $SS(f^{-1}F_i)\cap \pi^{-1}(W')=SS(G_i')\cap
\pi^{-1}(W')$, for all $i\in I$. On the other hand, one has a
distinguished triangle $$G'\to G\otimes\K_{H_\varepsilon}\to Rq_{1
!! }(K\otimes q_2^{-1}G_{H_\varepsilon})\xrightarrow{+1},$$ where
$K$ denotes the complex $\K_{s^{-1}(\gamma^a)}\to \K_{s^{-1}(0)}$,
with $\K_{s^{-1}(0)}$ in degree $0$. Moreover, one has
$$J(G')\simeq\underset{\underset{i\in I
}{\longrightarrow}}{\sideset{``}{"}\lim } J(G'_i),$$ $$J(Rq_{1 !!
}(K\otimes
q_2^{-1}G_{H_\varepsilon}))\simeq\underset{\underset{i\in I
}{\longrightarrow}}{\sideset{``}{"}\lim } J(Rq_{1 ! }(K\otimes
q_2^{-1}(G_i)_{H_\varepsilon})),$$ and $SS(Rq_{1 !}(K\otimes
q_2^{-1}(G_i)_{H_\varepsilon}))\cap
(Y\times\text{Int}(\gamma^{\circ a}))=\emptyset$, for all $i\in
I$. It follows that $G\simeq f^{-1}F$ in
$D^b(\text{I}(\K_Y;W'\times \text{Int}(\gamma^{\circ a}))$ and
$SS(f^{-1}F_i)\cap \pi^{-1}(W')=SS(G_i)\cap \pi^{-1}(W')$, for all
$i\in I$.

Let us now suppose dim$Z=m\geq 1$. We may assume $Z=\R^m$,
$X=\R^d$, $p=(0;\xi_0)$ and $\xi_0=(0,\xi_0')$. Set $X'=\R\times
X$ and let $h$ (resp. $g$) denote the projection $Y\to
X';(z',z'',x)\mapsto (z'',x)$ (resp. $X'\to X;(z'',x)\mapsto x$).
Since $Y\times_X T^*X\subset Y\times_{X'}T^*X'$, by induction,
there exist an open neighborhood $U'\subset U$ of $p$ and $F\in
D^b(\text{I}(\K_{X'}))$ such that $h^{-1}F\simeq G$ in
$D^b(\text{I}(\K_Y; U'))$ and $J(F)\simeq\underset{\underset{i\in
I }{\longrightarrow}}{\sideset{``}{"}\lim }J(F_i)$, for some
$F_i\in D^b(\K_{X'})$ with  $SS(h^{-1}F_i)\cap U'=SS(G_i)\cap U'$,
for all $i\in I$. Since, for each $i\in I$, $SS(h^{-1}F_i)=h_d
h_\pi^{-1}(SS(F_i))$, it follows that $SS(F_i)\cap h_\pi
h_d^{-1}(U')\subset X'\times_X T^*X$.

By the case $m=1$, we conclude that there exist an open
neighborhood $V\subset h_\pi h_d^{-1}(U')$ of $(0;0,\xi'_0)\in
T^*X'$ and $H\in D^b(\text{I}(\K_{X}))$ such that $g^{-1}H\simeq
F$ in $D^b(\text{I}(\K_{X'}; V))$ and
$J(H)\simeq\underset{\underset{i\in I
}{\longrightarrow}}{\sideset{``}{"}\lim }J(H_i)$, for some $H_i\in
D^b(\K_{X'})$ with $SS(g^{-1}H_i)\cap V=SS(F_i)\cap V$, for all
$i\in I$.

Set $V'=\{(z',z'',x;\xi',\xi'',\eta)\in T^*Y;
(z'',x;\xi'',\eta)\in V\}$. By Proposition \ref{SP:6}, one has
$f^{-1}H\simeq h^{-1}F$ in $D^b(\text{I}(\K_Y, V'))$ and
$SS(f^{-1}H_i)\cap V'=SS(h^{-1}g^{-1}H_i)\cap
V'=h_dh_\pi^{-1}(SS(g^{-1}H_i)\cap V)=h_dh_\pi^{-1}(SS(F_i)\cap
V)=SS(h^{-1}F_i)\cap V'$. The result follows.
\end{proof}

\subsection{Application to $\shd$-modules}

\hspace*{\parindent}Let $X$ be a complex manifold and let $\shm$
be a coherent $\shd$-module. We denote by $Char(\shm)$ the
characteristic variety of $\shm$, by $\sho_X^t$ the ind-sheaf of
temperate holomorphic functions on $X$ and we set for short:
$$Sol^t(\shm)=R\mathcal{I}{hom}_{\beta_X(\shd_X)}(\beta_X(\shm),
\sho_X^t).$$

It is proved in \cite{KS2} the equality
$$SS(Sol^t(\shm))=Char(\shm),$$ and that $Sol^t(\shm)$ is regular
if $\shm$ is regular holonomic. We shall now prove the following:

\begin{proposition}\label{Eg:1}
If $\shm$ is regular along an involutive vector subbundle $V$ of
$T^*X$, then $Sol^t(\shm)$ is regular along $V$.
\end{proposition}

\begin{proof}[\textbf{Proof}]
We may assume $X=Z\times Y$, for some complex manifolds $Z$ and $
Y$, that $f$ is the projection $X\to Y$ and that $V=X\times_Y
T^*Y$. By Lemma 3.6 of \cite{KMS} we may also assume
$\shm=\shd_{X\to Y}$. Then, by Lemma 7.4.8 of \cite{KS1}, one has:
$$Sol^t(\shm)\simeq f^{-1}\sho_Y^t,$$ and the result follows by
Proposition \ref{SP:9}.
\end{proof}

As another application we obtain a new proof of the following
result of \cite{KMS}:

\begin{theorem}
Let $V$ be an involutive vector subbundle of $T^*X$. Let $\shm$ be
a coherent $\shd_X$-module regular along $V$ and let $F\in
D^b_{\R-c}(\C_X)$. We have the estimate:
$$SS(R\mathcal{H}{om}_{\shd_X}(\shm, {thom}(F,\sho_X))\subset V
\widehat{+}SS(F)^a.$$
\end{theorem}

\begin{proof}[\textbf{Proof}]
We may assume $X=Z\times Y$, for some complex manifolds $Z$ and
$Y$, that $f$ is the projection $X\to Y$ and that $V=X\times_Y
T^*Y$.

One has $$R\mathcal{H}{om}_{\shd_X}(\shm, {thom}(F,\sho_X)\simeq
R\mathcal{H}{om}(F, Sol^t(\shm))\simeq$$ $$\simeq
\alpha_X(R\mathcal{I}{hom}(F, Sol^t(\shm))).$$

Arguing as in the proof of Proposition \ref{Eg:1}, it is enough to
prove the estimate
$$SS(R\mathcal{I}{hom}(F,f^{-1}\sho_Y^t))\subset V
\widehat{+}SS(F)^a.$$

Let us consider a filtrant inductive system $G_i$ in
$D^b(\text{Mod}^c(\K_X))$ such that $J(\sho_Y^t)\simeq
\underset{\underset{i\in I}{\longrightarrow}}{\sideset{``}{"}\lim
}J(G_i)$. Then:
$$J(R\mathcal{I}{hom}(F,f^{-1}\sho_Y^t))\simeq\underset{\underset{i}{\longrightarrow}}{\sideset{``}{"}\lim
}J(R\mathcal{I}{hom}(F,f^{-1}G_i))$$ and
$$SS(R\mathcal{I}{hom}(F,f^{-1}\sho_Y^t))\subset \bigcap_{J\subset
I}\overline{\bigcup_{j\in J}(SS(F)^a\widehat{+} SS(f^{-1}G_j))},$$
where $J$ runs over the family of cofinal subcategories of $I$
(see \cite{KS2}). Since $SS(f^{-1}G_j)\subset V$, for all $j\in
I$, the desired estimate follows.
\end{proof}

\subsection{Microsupport for sheaves on the subanalytic site}

\hspace*{\parindent}Let $\gamma$ be a closed convex proper cone in
$\R^n$ and $X\subset \R^n$ be an open subset. Denote by $q_1,
q_2:X\times X\to X$ the first and second projections and denote by
$s:X\times X\to X$ the map $(x,y)\mapsto x-y$. For an open
subanalytic subset $W\subset X$ we  define the functor
$\Phi_{\gamma, W}$ in $D^b(\K_{X_{sa}})$ by setting: $$
\Phi_{\gamma, W}(F)=Rq_{1 !!}(\rho_*\K_{s^{-1}(\gamma)\cap
q_1^{-1}(W)\cap q_2^{-1}(W)}\otimes q_2^{-1}F),$$ and we write
$\Phi_\gamma$ instead of $\Phi_{\gamma, X}$.

\begin{lemma}\label{SL:2}
Let $X$ be an open subset of $\R^n$ and let $F\in
D^b(\K_{X_{sa}})$. For every closed convex proper cone $\gamma$ in
$\R^n$, one has $$\Phi_{\gamma}(\I(F))\simeq
\I(\Phi_{\gamma}(F)).$$
\end{lemma}

\begin{proof}[\textbf{Proof}]
Recall that, for any real analytic map $f:X\to Y$, the functor
$\I$ commutes with the functors $Rf_{!!}$ and $f^{-1}$. Hence,
given a closed convex proper cone $\gamma$ in $\R^n$, one has:
$$\Phi_{\gamma}(\I(F))=Rq_{1 !!}(\K_{s^{-1}(\gamma)}\otimes
q_2^{-1}\I(F))\simeq  Rq_{1 !!}(\I(\K_{s^{-1}(\gamma)}\otimes
q_2^{-1}F))\simeq$$ $$\simeq \I(Rq_{1
!!}(\K_{s^{-1}(\gamma)}\otimes q_2^{-1}F))=\I(\Phi_{\gamma}(F)).$$
\end{proof}

Let $F\in D^b(\K_{X_{sa}})$. The microsupport of $F$, denoted by
$SS(F)$, is the closed conic subset of $T^*X$ defined by the
equality: $$SS(F)=SS(\I(F)).$$

\begin{proposition}\label{SP:11}
Assume $X$ is an open subset of $\R^n$. Let $F\in
D^b(\K_{X_{sa}})$ and $p\in T^*X$, with $p=(x_0;\xi_0)$. The
conditions $(i)-(iv)$ bellow are all equivalent.

(i) There exists a conic open neighborhood $U$ of $p$ in $T^*X$
such that for any $G\in D^b_{\R-c}(\K_X)$, with
$\mathrm{supp}(G)\subset\subset \pi(U)$ and $SS(G)\subset U\cup
T_X^*X$, one has $\mathrm{Hom}_{D^b(\K_{X_{sa}})}(R\rho_*G,F)=0$.

(ii) There exist a closed convex proper cone $\gamma$ in $\R^n$,
with $\xi_0\in \mathrm{Int}(\gamma^{\circ})$, and a relatively
compact open subanalytic neighborhood $W$ of $x_0$, such that
$\Phi_{\gamma, W}(F)\simeq 0.$

(iii) There exist $F'\in D^b(\K_{X_{sa}})$ with compact support,
$F'\simeq F$ in an open subanalytic neighborhood of $x_0$ and a
closed convex proper cone $\gamma$ as in $(ii)$ such that
$\Phi_{\gamma}(F')\simeq 0$ in a neighborhood of $x_0$.

(iv) There exist $F'\in D^b(\K_{X_{sa}})$ with compact support
such that $F'\simeq F$ in an open subanalytic neighborhood of
$x_0$ and a closed convex proper cone $\gamma$ as in $(ii)$ such
that $Rq_{2 *}R\Gamma_{s^{-1}(\gamma)}(q_1^!F')\simeq 0$ in a
neighborhood of $x_0$.
\end{proposition}

\begin{proof}[\textbf{Proof}]
The proof of \noindent$(i) \Leftrightarrow (iv)$ is analogous to
the proof of Lemma \ref{SP:2}.

\noindent$(i)\Leftrightarrow (ii)$ follows by Lemma 4.1 of
\cite{KS2} together with Proposition \ref{SP:2}. In fact, for each
$G\in D^b_{\R-c}(\K_X)$ with compact support, we may choose $G'\in
D^b(\text{Mod}_{\R-c}^c(\K_X))$ quasi-isomorphic to $G$. From the
equivalence of categories $$D^b(\K_{X_{sa}})\simeq
D^b(\text{I}\R-c(\K_X))\xrightarrow{\I}
D^b_{\text{I}{\R-c}}(\text{I}(\K_X)),$$ and since
$D^b_{\text{I}\R-c}(\text{I}(\K_X))$ is a full subcategory of
$D^b(\text{I}(\K_X))$, one gets that:
$$\text{Hom}_{D^b(\K_{X_{sa}})}(R\rho_*G,F)\simeq\text{Hom}_{D^b(\K_{X_{sa}})}(\rho_*G',F)\simeq$$
$$\simeq \text{Hom}_{D^b(\text{I}\R-c(\K_X))}(\rho_*
G',F)\simeq\text{Hom}_{D^b_{\text{I}{\R-c}}(\text{I}(\K_X))}(\I(\rho_*G'),\I(F))\simeq$$
$$\simeq \text{Hom}_{D^b(\text{I}(\K_X))}(G',\I(F))\simeq
\text{Hom}_{D^b(\text{I}(\K_X))}(G,\I(F)).$$

\noindent$(ii) \Rightarrow (iii)$ is obvious by taking $F_W$ as
$F'$.

\noindent$(iii) \Rightarrow (ii)$ follows by Lemma 4.1 of
\cite{KS2} together with Lemma \ref{SL:2}.
\end{proof}

\begin{corollary}\label{SC:1}
Let $F\in D^b(\K_{X_{sa}})$. The microsupport of $F$ is the closed
conic subset of $T^*X$ complementary of the set of points $p\in
T^*X$ such that one of the equivalent conditions of Proposition
\ref{SP:11} is satisfied.
\end{corollary}

\begin{proof}[\textbf{Proof}]
It is a consequence of Lemma 4.1 of \cite{KS2} and Proposition
\ref{SP:11}.
\end{proof}

\begin{proposition}\label{SP:12}
(i) For $F\in D^b(\K_{X_{sa}})$, one has $SS(F)\cap
T_X^*X=\mathrm{supp}(F)$.

(ii) Let $F\in D^b(\K_{X_{sa}})$, then $SS(\rho^{-1}F)\subset
SS(F)$.

(iii) Let $F\in D^b(\K_X)$, then $SS(F)\subset SS(R\rho_*F)$.

(iv) Let $F_1\to F_2\to F_3\xrightarrow{+1}$ be a distinguished
triangle in $D^b(\K_{X_{sa}})$. Then for $i,j,k\in\{1,2,3\}$,
$$SS(F_i)\subset SS(F_j)\cup SS(F_k), \ \text{for $j\neq k$}.$$
\end{proposition}

\begin{proof}[\textbf{Proof}]
\noindent$(i)$ is a consequence of Proposition 4.3 of \cite{KS2}.

\noindent$(ii)$  One has $\rho^{-1}\simeq \alpha\circ \I$ and by
Proposition 4.3 of \cite{KS2}:
$$SS(\rho^{-1}F)=SS(\alpha(\I(F))\subset SS(\I(F))=SS(F).$$

\noindent$(iii)$ follows from $(iii)$:
$$SS(F)=SS(\rho^{-1}R\rho_*F)\subset SS(R\rho_*F).$$

\noindent$(iv)$ follows from Proposition 4.3 of \cite{KS2}.
\end{proof}

\begin{proposition}\label{SP:13}
Let $F\in D^b(\K_{X_{sa}})$ and $G\in D^b(\K_{Y_{sa}})$. Then:
$$SS(F\boxtimes G)= SS(F)\times SS(G).$$
\end{proposition}

\begin{proof}[\textbf{Proof.}]
Since $\I(F\boxtimes G)\simeq \I(F)\boxtimes \I(G)$, it follows by
Proposition \ref{SP:4}. In fact, one has: $$SS(F\boxtimes
G)=SS(\I(F\boxtimes G))=SS(\I(F))\times SS(\I(G))=SS(F)\times
SS(G).$$
\end{proof}

\begin{proposition}\label{SP:15}
Let $M$ be a closed submanifold of $X$ and let $j$ denote the
closed embedding $M\hookrightarrow X$. Let $G\in
D^{b}(\mathbf{k}_{M_{sa}})$. Then, $$SS(Rj_*G)=
j_{\pi}j_d^{-1}(SS(G)).$$
\end{proposition}

\begin{proof}[\textbf{Proof}]
Since $Rj_*\simeq Rj_{!!}$ and $\I$ commutes with $Rj_{!!}$, the
result follows by Proposition \ref{SP:5}.
\end{proof}

\begin{proposition}\label{SP:16}
Let $Y$, $X$ be real analytic manifolds, $f\colon Y\to X$ be a
smooth morphism and $F\in D^b (\K_{X_{sa}})$. Then $$SS(f^{-1}F)=
f_d f_{\pi}^{-1}(SS(F).$$
\end{proposition}

\begin{proof}[\textbf{Proof}]
It follows by Proposition \ref{SP:13}.
\end{proof}

\begin{proposition}\label{SP:17}
Let $F\in D^b(\K_{X_{sa}})$ and $G\in D^b(\K_{Y_{sa}})$. Then:
$$SS(R\mathcal{H}{om}(q_2^{-1}G,q_1^{!}F))\subset SS(F)\times
SS(G)^a.$$
\end{proposition}

\begin{proof}[\textbf{Proof}]
The proof is an adaptation of the proof of Proposition \ref{SP:7},
using Propositions 5.4.6, 5.4.7, 5.4.5 and 5.2.9 of \cite{L}
instead of Propositions 5.3.8, 5.3.10, 5.3.5 and 5.2.9 of
\cite{KS1}, respectively.
\end{proof}

{\small Ana Rita Martins

Centro de {\'A}lgebra da Universidade de Lisboa, Complexo 2,

Avenida Prof. Gama Pinto,

1699 Lisboa Portugal

arita@mat.fc.ul.pt}
\end{document}